\documentclass[a4paper,12pt]{amsart}
\usepackage{amsmath,amsfonts,amssymb}
\usepackage{amsthm}
\usepackage{graphics,graphicx,epsf}
\usepackage[usenames]{color}
\usepackage{color}
\usepackage{hyperref}
\usepackage{refcheck}  
\usepackage{tikz}

\newtheorem{theorem}{Theorem}[section]
\newtheorem{definition}[theorem]{Definition}
\newtheorem{proposition}[theorem]{Proposition}
\newtheorem{lemma}[theorem]{Lemma}
\newtheorem{corollary}[theorem]{Corollary}
\newtheorem{remark}[theorem]{Remark}

\renewcommand{\proof}{{\noindent \bf Proof. \ }}

\newcommand{\R}{\mathbb{R}}
\newcommand{\Sprocess}{ \{S(t,\tau); t \geq \tau\in\mathbb{R} \}}
\newcommand{\Sautprocess}{\{S_0(t,\tau); \ t \geq \tau\in\mathbb{R} \}}
\newcommand{\Sdyn}{\{S_0(t)\ ; \ t \geq 0 \}}
\newcommand{\sgn}{\textrm{sgn}}

\numberwithin{equation}{section}

\title[nonlocal problems]{Asymptotic behavior for a class of nonlocal nonautonomous problems}

\author[F. D. M. Bezerra]{Flank  D. M. Bezerra$^1$}
\address[F. D. M. Bezerra]{Departamento de Matem\'atica, Universidade Federal da Para\'{\i}ba,
Cidade Universit\'{a}ria-Campus I, 58051-900 Jo\~{a}o Pessoa PB, Brazil.}
\email{flank@mat.ufpb.br}

\author[S. H. Da Silva]{Severino H. da Silva$^2$}\thanks{$^2$Research partially
supported by CAPES/CNPq}
\address[S. H. Da Silva]{Unidade Acad\^emica de Matem\'atica,
Universidade Federal de Campina Grande, 58051-900 Campina Grande PB,
Brazil.} \email{horacio@dme.ufcg.edu.br}

\author[A. L. Pereira]{Ant\^onio L. Pereira$^3$}
\address[A. L. Pereira]{Universidade de S\~{a}o
Paulo, Instituto de Matem\'atica e Estat\'{\i}stica, 13565-905 S\~{a}o
Paulo SP, Brazil.}
\email{alpereir@ime.usp.br}

\date{\today}

\begin{document}

 \maketitle

\begin{abstract}
In this paper we consider the nonlocal nonautonomous evolution problem
\[
\begin{cases}
\partial_t u    =- u  + g \left(t,  Ku  \right)   \ \ in \ \  \Omega,\\
u  = 0  \ \  in \ \
\mathbb{R}^N\backslash\Omega.
\end{cases}
\]
where $\Omega$ is a smooth bounded domain in $\mathbb{R}^N$, $g: \mathbb{R} \times \mathbb{R} \to \mathbb{R}$ and $K$ is an integral operator with a symmetric kernel. We prove  existence and some regularity properties of the  pullback attractor. We also  show additional forward  asymptotic results in the asymptotically autonomous case, using the properties of the Lyapunov functional for the limiting problem.

\end{abstract}



\section{Introduction}

The understanding of the dynamics of nonautonomous evolution equations has attracted the attention of many researchers in recent years; see for instance \cite{CLR,CS,CV,Kl,Sch} and \cite{Sell}.
   The results obtained here constitute an extension of the ones in
 \cite{FAH} to our (more general) situation. We  also
 prove new results in a different phase space ($L^{\infty}$)  and in
the asymptotically autonomous case.

More precisely, we consider here the nonlocal nonautonomous evolution problem
\begin{equation}\label{prob}
\begin{cases}
\partial_t u (t,x) =- u (t,x)  + g(t, Ku(t,x))  \ \ \textrm {for }
 \ \ t \geq  \tau \in \R \ \ \textrm{and}    \ \ x \in  \Omega, \\
u(\tau, \cdot) = u_\tau (\cdot)   \ \  \textrm{in}  \ \ \Omega, \\
u(t, x)   = 0  \ \  \textrm{for}   \ \ t \geq  \tau \in \R
 \ \ \textrm{and}  \ \  x \in
\mathbb{R}^N\backslash\Omega.
\end{cases}
\end{equation}
where $\Omega$ is a bounded smooth domain in $\mathbb{R}^N$ ($N\geq1$),
$g: \mathbb{R} \times \mathbb{R} \to \mathbb{R}$ is a (sufficiently smooth)   function and $K$ is an integral operator with symmetric kernel
\begin{equation*}
Ku(\cdot,x):=\int_{{\mathbb{R}^{N}}} J(x,y)u(\cdot,y)dy.
\end{equation*}

 We will suppose, without loss of generality,  that $\int_{\mathbb{R}^N} J(x,y) \, d \, y
 = \int_{\R^N} J(x,y) \, d \, x = 1.$
\medskip

\medskip

The paper is organized as follows. After recalling some concepts and results about attractors for infinite-dimensional nonautonomous dynamical
systems in Section \ref{prelim}, we prove  well posedness of \eqref{prob} in a Banach space which is isomorphic to $L^p(\Omega)$, $1\leq p\leq\infty$, in Section \ref{wellposed}, and existence and some regularity properties for the pullback attractor on Section \ref{PullAttractors}. In Section \ref{comparison} we prove some comparison results. In Section \ref{Sec6} is dedicated to the proof of upper semi-continuity of the family of pullback attractors with respect to the functional parameter $g$. Finally, on Section \ref{NS}, we prove some additional results on the asymptotic behavior of \eqref{prob} for some special cases,  including the asymptotic autonomous one.

\section{Some preliminary definitions and results} \label{prelim}

We start by recalling  the concepts of  evolution processes  and pullback attractors which proved to be a useful tool to study the asymptotic dynamics of infinite dimensional nonautonomous dynamical system. We refer to \cite{Carvalho,CLR,CS,CV,Kl} and \cite{Sch}  for more details and proofs.

\begin{definition}[Nonlinear evolution process]
A {\rm nonlinear evolution process} in a Banach space $\mathcal{X}$ is a family of
  maps $\{S(t,\tau); t\geq\tau\in\mathbb{R}\}$ (not necessarily linear) from $\mathcal{X}$ into itself with the following properties:
\begin{itemize}
\item[(1)] $S(t,t)=I$, for all $t\in\mathbb{R}$,
\item[(2)] $S(t,\tau)=S(t,s)S(s,\tau)$, for all $\tau\leq s\leq t$,
\item[(3)] The map $\{(t,\tau)\in\mathbb{R}^2;\ t\geq\tau\}\times \mathcal{X}\ni (t,\tau,x)\mapsto S(t,\tau)x\in \mathcal{X}$ is continuous.
\end{itemize}
\end{definition}

In the particular case where each $S(t,\tau)$ is linear, $t\geq\tau\in\mathbb{R}$, we say that $\{S(t,\tau); t\geq\tau\in\mathbb{R}\}$ is a linear process.

If $S(t,\tau)={S}(t-\tau,0)$ for all $t\geq\tau\in\mathbb{R}$, we say that the process is \emph{autonomous}, and the family
  $\{\bar{S}(t)= S(t,0); t\geq0\}$ is then called a \emph{semigroup} or a \emph{dynamical system}.

 \begin{definition}
 The pullback (resp. forward) orbit of a subset $B$ of $\mathcal{X}$, at time $t \in \R$, is the set
 $\gamma_p(B,t) : = \cup_{s \leq t} S(t,s) B $  (resp.  $\gamma_f(B,t) : = \cup_{s \geq t} S(s,t) B $).
\end{definition}

\begin{definition}
A {\rm globally-defined solution}  (or simply a {\rm global solution}) of the  evolution process $\{S(t,\tau); t\geq\tau\in\mathbb{R}\}$ through $\psi_0\in \mathcal{X}$ at
 time $ \tau_0$ is a function $\psi:\mathbb{R}\to \mathcal{X}$ such that
 $\psi(\tau_0)=\psi_0$,  and  $S(t,\tau)\psi(\tau)=\psi(t)$,
 for all $t\geq\tau$.
 A global solution through  $\psi_0\in \mathcal{X}$ is a
  global solution through  $\psi_0\in \mathcal{X}$ at some time.
\end{definition}

\begin{definition}[Pullback Attraction]
A family of sets $\{ \mathcal{K}(t); t\in\mathbb{R}\}$ {\rm pullback attracts bounded subsets of}  $\mathcal{X}$ under $\{S(t,\tau); t\geq\tau\in\mathbb{R}\}$, at time
 $t$  if,
for each bounded set $C\subset\mathcal{X}$,
\[
\lim_{\tau\to-\infty}\mathrm{dist}(S(t,\tau)C,\mathcal{K}(t))=0,
\]
where $\mathrm{dist}(\cdot,\cdot)$ denotes the Hausdorff semi-distance in $\mathcal{X}$,
\begin{equation}\label{DD}
\mathrm{dist}(A,B)=\sup_{a\in A}\inf_{b\in B}|a-b|_{\mathcal{X}}.
\end{equation}
We observe that the Hausdorff semi-distance between $A$ and $B$, $\mathrm{dist}(A, B)$, examines
how the set $A$ is contained in the set $B$. For example, $\mathrm{dist}(A, B)=0$ if and only if
$A$ is contained in the closure of the set $B$.

A family of sets $\{\mathcal{K}(t); t\in\mathbb{R}\}$ {\rm pullback attracts bounded subsets of}  $\mathcal{X}$ under $\{S(t,\tau); t\geq\tau\in\mathbb{R}\}$ if
 it   { pullback attracts bounded subsets of}  $\mathcal{X}$ under $\{S(t,\tau); t\geq\tau\in\mathbb{R}\}$, at time
 $t$, for any $t$.

\end{definition}

 \begin{definition}
A family of bounded subsets $\{B(t); t\in\mathbb{R}\}$ of $\mathcal{X}$   {\rm
  pullback absorbs bounded subsets of }  $\mathcal{X}$ under $\{S(t,\tau); t\geq\tau\in\mathbb{R}\}$ at time $t \in \R.$
 if, for any bounded set $C\subset \mathcal{X}$ and $\tau\leq t$, there exists $\tau_0(\tau, C)\in\mathbb{R}$ such that
$$
S(\tau,\tau-r)C\subset B(t),\ \mbox{for all}\ r\geq\tau_0 (\tau, C).
$$
A family of bounded subsets $\{B(t); t\in\mathbb{R}\}$ of $\mathcal{X}$   {\rm
  pullback absorbs bounded subsets of }  $\mathcal{X}$ under $\{S(t,\tau); t\geq\tau\in\mathbb{R}\}$ if it    {
  pullback absorbs bounded subsets of }  $\mathcal{X}$ under $\{S(t,\tau); t\geq\tau\in\mathbb{R}\}$ at time $t$, for any $t \in \R$.
\end{definition}
\begin{remark}
It is clear that, if  $\{B(t); t\in\mathbb{R}\}$    {
  pullback absorbs bounded subset of }  $\mathcal{X}$ at time $t$, then
    it  {
  pullback attracts  bounded subsets of }  $\mathcal{X}$ at time $t$.

When $B(t) = B,  \ \mbox{ for all} \in \R  \ $, where $B$ is
 some fixed bounded set, it is also said  that
 the set $B$ (instead of the family  $\{B(t); t\in\mathbb{R}\}$)  pullback
 absorbs bounded  subsets  of  $\mathcal{X}$.
\end{remark}

\begin{definition}
An evolution process $\{S(t,\tau); t\geq\tau\in\mathbb{R}\}$ is said to be {\rm strongly pullback  bounded dissipative} if there exists a family
of sets $\{\mathcal{K}(t); t\in\mathbb{R}\}$  that  pullback attracts bounded subsets of  $\mathcal{X}$ under $\{S(t,\tau)\}$ at time $s$ for any $s\leq t$.
\end{definition}

Next, we present the notion of a  set  of pullback asymptotic states.

\begin{definition}[Pullback attractor] \label{def-pullback-attract}
A family $\{\mathcal{A}(t); t\in\mathbb{R}\}$ of compact subsets of $\mathcal{X}$ is said to be a {\rm pullback attractor} for the evolution process $\{S(t,\tau); t\geq\tau\in\mathbb{R}\}$ if it is invariant, i.e., $S(t,\tau)\mathcal{A} (\tau) = \mathcal{A} (t )$ for all $\tau\leq t$, pullback attracts bounded subsets of $\mathcal{X}$, and is minimal, that is, if there is another family of closed sets $\{C(t); t\in\mathbb{R}\}$ which pullback attracts bounded subsets of $\mathcal{X}$, then $\mathcal{A}(t) \subset C(t)$, for all $t\in\mathbb{R}$.
\end{definition}

\begin{remark}
The minimality requirement in the Definition \ref{def-pullback-attract} is an
 addition with respect to the theory of attractors for semigroups and  is necessary to ensure uniqueness (see \cite{CLR}). It can be dropped if we require that $\bigcup_{\tau \leq t} \mathcal{A} (\tau )$ is bounded for any $t \in \R$. In this case, we also have
$$
\mathcal{A}(t) =\Big\{  \xi(t) : \xi: \R \to \mathcal{X}\ \textrm{is a global backwards bounded solution of } S(t,\tau)  \Big\},
$$
for all $t\in\mathbb{R}$.
\end{remark}

\begin{definition}
An evolution process $\{S(t,\tau); t\geq\tau\in\mathbb{R}\}$ in a Banach space $\mathcal{X}$ is  {\rm pullback asymptotically compact} if, for each $t\in\mathbb{R}$, each sequence $\{\tau_k\}_{k\in\mathbb{N}}$ in $(-\infty, t ]$ such that $\tau_k\to-\infty$ as $k\to\infty$, and each bounded sequence $\{z_k\}_{k\in\mathbb{N}}$ in $\mathcal{X}$ with $\{S (t , \tau_k )z_k\}_{k\in\mathbb{N}}$ bounded, the sequence $\{S (t , \tau_k )z_k\}_{k\in\mathbb{N}}$ possesses a convergent subsequence.
\end{definition}


The following two results are proved in \cite{CLR} (see also \cite{CS}) and will be used to prove the existence of the pullback attractor for the evolution process generated by (\ref{prob}), (see Section \ref{PullAttractors}).

\begin{theorem} \label{critcompact}
If an evolution process $\{S(t,\tau); t\geq\tau\in\mathbb{R}\}$ in a Banach space $\mathcal{X}$  satisfies the properties
\[
S(t,\tau)=T(t,\tau)+U(t,\tau),\ t\geq\tau,
\]
where $U(t,\tau)$ is a compact operator and there exists a non-increasing function $k:[0,+\infty)\times[0,+\infty)\to\mathbb{R}$ with $k(\sigma,r)\to0$ as $\sigma\to+\infty$, and for all $\tau\leq t$ and $z\in \mathcal{X}$ with $|z|_{\mathcal{X}}\leq r$, $|T(t,\tau)|_{\mathcal{X}}\leq k(t-\tau,r)$, then the process $\{S(t,\tau); t\geq\tau\in\mathbb{R}\}$ is pullback asymptotically compact.
\end{theorem}

\begin{theorem}\label{Teorema 0.2.1}
If an evolution process $\{S(t,\tau); t\geq\tau\in\mathbb{R}\}$ in a Banach space $\mathcal{X}$ is strongly pullback  bounded dissipative and pullback asymptotically compact, then $\{S(t,\tau); t\geq\tau\in\mathbb{R}\}$ possesses a pullback attractor $\{\mathcal{A}(t); t\in\mathbb{R}\}$. Moreover, the union $\bigcup_{\tau\leq t}\mathcal{A}(\tau)$ is bounded for each $t\in\mathbb{R}$, and the global pullback attractor is given by
\[
\mathcal{A}(t)=\Omega_\wp(\mathcal{B}(t),t):=\bigcap_{s\leq
t}\overline{\bigcup_{\tau\leq s}S(t,\tau)\mathcal{B}(t)},\ \mbox{for
each}\ t\in\mathbb{R},
\]
where $\{B(t); t\in\mathbb{R}\}$ is a family of bounded subsets of $\mathcal{X}$
which pullback attracts bounded subsets of $\mathcal{X}$ at time $\tau$, for any
 $\tau \leq t$.
\end{theorem}

This theorem extends the analogous result for semigroups (cf. Theorem 1.1, Chapter 1 in \cite{T}). The pullback attractor of strongly bounded dissipative process however, is always bounded in the past. To be more precise, for every $t\in\mathbb{R}$ the union $\bigcup_{\tau\leq t}\mathcal{A}(\tau)$ is bounded in $\mathcal{X}$.

\section{Well posedness}\label{wellposed}

 In this section   we  show  that the problem (\ref{prob}) is well posed in
  suitable phase spaces.



 \vspace{3mm}

 Consider for any $1 \leq p \leq \infty$  the  subspace  $X$ of $L^{p}(\mathbb{R}^{N})$   given  by
 \begin{equation*}
 X= \left\{ u \in  L^{p}(\mathbb{R}^{N});\ \ u(x)= 0, \textrm{ if }
 x\in \mathbb{R}^N\backslash\Omega \right\}
 \end{equation*}
with the induced norm.  The space $X$ is canonically isomorphic to
 $L^p(\Omega)$ and we usually identify the two spaces, without further comment.  We also use the same notation for a
 function in   $ \mathbb{R}^N $ and its restriction to  $\Omega$    for simplicity, wherever we believe the intention is clear from the context.


In order to obtain well posedness of (\ref{prob}), we  consider the
 Cauchy problem
\begin{equation}\label{CP}
\begin{cases}
\displaystyle { \partial_t u  } = -u +  F(t,u), \\
 u(\tau)=u_{\tau},
\end{cases}
\end{equation}
where the map $F: \R \times X  \to  X $
  is defined by
\begin{equation} \label{mapF}
F(t,u)(x)=
\begin{cases}
 g(t, Ku(t,x)), & x\in \Omega,\\
 0, & x\in\mathbb{R}^N\backslash\Omega.
 \end{cases}
\end{equation}

\begin{definition}

 A {\rm solution} of  (\ref{CP}) in $[\tau,s), \ \tau < s$  is a  continuous function $u:[\tau,s) \to X  $   such that $u(\tau)= u_{\tau} $,
 the derivative with respect to $t$ exists and $\partial_t u(t,\cdot)$ belongs to $X$,
 and the differential equation in
  (\ref{CP}) is satisfied for $ t \in [\tau, s)$.
\end{definition}

  The  map  given by
  \begin{equation} \label{mapK}
 K u (\cdot,x)  :=   \int_{\R^N} J(x,y)  u(\cdot,y)d y
\end{equation}
 is well defined as a bounded linear operator in  various function spaces, depending on the properties assumed for $J$.  We collect here some estimates for this map
 which  will be  used in the sequel.

\begin{lemma} \label{boundK}
  Let $K$ be the map defined  by (\ref{mapK}) and    $\|J\|_{p}$:=
 $\sup_{ x \in \Omega } \|J(x,\cdot)\|_{L^{p} (\Omega) }$, $ 1 \leq p \leq \infty.$
  If $u \in   L^p{(\Omega)}, \ 1 \leq p \leq \infty$, then
 $Ku \in L^p{(\Omega)}$, and

   \begin{equation}  \label{estimateL1}
   \|Ku\|_{L^{p}(\Omega)}   \leq   \|J\|_{1}
   \| u\|_{L^{p}(\Omega)} =  \| u\|_{L^{p}(\Omega)}.
 \end{equation}
   \begin{equation}  \label{estimateLq}
 |Ku (x)|  \leq   \|J\|_{q}      \| u\|_{L^p{(\Omega)}},
  \mathrm{for\ all}\  x \in \Omega,
 \end{equation}
\textrm{ where $1\leq q \leq \infty$
 is the conjugate exponent of $p$}.
\end{lemma}

\proof
Estimate  (\ref{estimateL1}) follows  from  Young's inequality and the fact that
   $\|J\|_{1} = 1 $  and estimate
   (\ref{estimateLq}) follows from   H\"{o}lder inequality.

\qed

\begin{definition} \label{loclip}
If $E$ is a normed space, and $I \subset \R$ is an interval,
we say that a function  $F : I \times E \to E$ is
   \emph{locally Lipschitz
 continuous  (or simply locally Lipschitz) in the second variable } if,
 for any $ (t_0, x_0) \in  I \times E$, there exists  a constant
 $C$ and a rectangle $R = \{ (t,x) \in I \times E \ : \
 |t-t_0| <a,  \|x-x_0  \|< b \}$ such that,
 if $ (t,x)$ and $ (t, y)$ belong to $R$   then
 $ \|F(t,x) - F(t,y) \| \leq  C \| x-y \| $; we say that $F$ is
   \emph{ Lipschitz
 continuous on bounded sets  in the second variable } if the rectangle
 $R$ in the previous definition can chosen as any bounded rectangle in
$I \times E$
 \end{definition}



\begin{remark}
  The two definitions in \ref{loclip}  are equivalent if the normed space $ E  $ is locally compact.
\end{remark}


\begin{proposition}\label{WellP}
Suppose,  in addition to the hypotheses of Lemma \ref{boundK}, that the
function $g$ is  locally Lipschitz continuous in the second variable
 in $\R \times \R$.
Then
  the function $F$ defined by \eqref{mapF} is bounded
 Lipschitz continuous in the second variable in
  $\R \times L^{p}(\Omega), 1  \leq p \leq \infty$.
\end{proposition}

\proof
Suppose $(t_0, u_0) \in \R \times X$. Then if $R$ is the rectangle
  $ R: = \{ (t,u) \in \R \times X \ |  \ |t-t_0| < a , \
 \|u-u_0 \|_{L^{p}(\Omega)} < b $
  it follows  from (\ref{estimateLq})  that
$|  K u_0 (x) |  <    \|J\|_{q} \|u_0 \|_{L^{p}(\Omega)} $,
$| K u (x) - K u_0 (x) |  <    \|J\|_{q} b $,
for any $x \in \Omega.$
 Let
  $k_{R'} $ be  the Lipschitz constant of $g$ in the
 rectangle
 $ R':=  \{  (t, x) \in I \times \R \ | \  |t-t_0| < a , \
  \| x \| \leq    \|J\|_{q} (\|u_0 \|_{L^{p}(\Omega)} + b)  \}$.
 Then, if  $  ( t, u), (t,v) \in R $    we obtain
   $ |g(t, Ku(x) )  - g( t, Kv(x) ) | \leq k_R'   | Ku(x)  - Kv(x)  |   $
 for any $x \in \Omega.$
 For   $ 1 \leq p < \infty$, it follows then,
from (\ref{estimateL1}),  that
 \begin{eqnarray*}
 \|g(t, Ku)  - g( t, Kv )\|_{L^p(\Omega)}
  & \leq  &  \left\{
 \int_{\Omega}    k_{R' }  | Ku(x)  - Kv(x)  |^p \,  dx\right\}^{1/p} \\
 &  =  &
 k_{R'}  \left\{  \int_{\Omega}     | K (u-v)(x))  |^p \,  dx\right\}^{1/p}  \\
  &  =  &
 k_{R'}   \| K (u-v) \|_{L^p(\Omega)} \\
 &  =  &
 k_{R'}   \| u-v \|_{L^p(\Omega)}.
  \end{eqnarray*}
 If $p =\infty$ the same inequality follows immediately from
 (\ref{estimateL1}).
 Thus,  the map
\[u \in X \mapsto  g(t, K u) \in X \]  is  Lipschitz in the
 rectangle $R$  and, therefore, so is the map $F(t, \cdot)$.

\qed

\bigskip

From the result above, it follows from well known results that  the  problem  (\ref{CP}) has a local solution for any initial condition in $X$.
 For the  global existence, we need the following result
 (\cite{ladas} - Theorem 5.6.1)

\begin{theorem} \label{ladas} Let $X$ be a Banach space, and suppose that
 $g: [ t_0, \infty)  \times X \to X$  is continuous and
  $ \|g(t,u) \|  \leq h(t, \|u  \|);  \textrm{ for all } (t,u) \in
    [t_0, \infty) \times X $, where   $h: [t_0, \infty)
 \times  \R^{+}  \to   \R^{+} $ is continuous and
$ h(t,r)$  is  non decreasing in $r \geq 0$, for each
$t \in [ t_0, \infty) $. Then, if  the maximal solution
 $r(t,t_0,r_0) $ of the scalar initial value problem
 \[  r'= h(t,r), \quad r(t_0) = r_ 0, \]
 exists throughout $[t_0, \infty)$,
 the maximal interval of existence of any  solution
 $u(t,t_0,y_0)$ of the initial value problem
 \[ \frac{du}{dt} = g(t,u), \quad t \geq t_0, \quad u(t_0) = u_0,   \]
 also contains   $[ t_0, \infty)$.

\end{theorem}

  \begin{corollary} \label{globalexist}
   Suppose,  in addition to the hypotheses of Proposition \ref{WellP}, that
  $ g $  satisfies:
\begin{equation}\label{dissip1}
 \limsup_{|x| \to \infty} \frac{|g(t, x)|}{|x|} < k_1,
\textrm{ for some constant }
  k_1 \in \R.
 \end{equation}
Then the   problem  (\ref{CP})  has a unique   globally defined
 solution for any initial
 condition in  $X$,
 which is  given for $t\geq\tau$ by
 the ``variation of constants formula''

 \begin{eqnarray}\label{EP_1}
 u(t,\tau,x;u_{\tau}) & = &
 e^{-(t-\tau)}u_{\tau}(x)+
 \displaystyle\int_{\tau}^{t}e^{-(t-s)} F(s, u(s, \tau, x ; u_{\tau})) \, ds
    \quad  x\in \mathbb{R}^N \nonumber\\
   &  = &
 \begin{cases}
 e^{-(t-\tau)}u_{\tau} (x)+\displaystyle\int_{\tau}^{t}e^{-(t-s)}g(s,
K(u(s,\tau, \cdot;
 u_{\tau}))(x))\, ds, & x\in \Omega, \\
  0, & x\in \mathbb{R}^N \backslash\Omega.
  \end{cases}
 \end{eqnarray}
\end{corollary}
\proof
From Proposition \ref{WellP}, it follows that the right-hand-side of
 (\ref{CP}) is  Lipschitz continuous in bounded sets of
 $X$ and, therefore, the Cauchy
problem (\ref{CP}) is well posed in $ X $ with a unique local solution $u(t,\tau,x;u_\tau)$,
 given by (\ref{EP_1})  (see \cite{DK}).

From condition  (\ref{dissip1}), it follows that there is a constant $k_2$,
 such that
 \begin{equation} \label{dissip2}
 \| g(t,x)  \| \leq k_2 + k_1 |x |, \ \ \textrm{ for any } \ \ x \in \R.
 \end{equation}

 If $ 1 \leq p < \infty$,  we obtain from   (\ref{dissip2})  and
  (\ref{estimateL1})  that
 \begin{eqnarray*}
  \| g(t,  K u )\|_{L^p(\Omega)} &   \leq &   k_2 |\Omega|^{1/p}    +
  k_1 \|   K u \|_{L^p(\Omega)} \\
 & \leq &  k_2 |\Omega|^{1/p}    +
  k_1 \|    u  \|_{L^p(\Omega)}.
\end{eqnarray*}
 For $p= \infty$, we obtain by the  same arguments (or by making
 $p \to \infty$), that
 \begin{eqnarray*}
  \| g(t,  K u ) \|_{L^{\infty}(\Omega)} &\leq  k_2     +
  k_1 \|    u  \|_{L^{\infty}(\Omega)}.
\end{eqnarray*}

  Defining $ h: [t_0, \infty)
 \times  \R^{+}  \to   \R^{+} $, by $h(t,r) =   k_2  |\Omega|^{1/p}   +
 (k_1 +1) r $, it follows that Problem \ref{CP}
  satisfies the hypothesis of  Theorem \ref{ladas}  and the global existence
 follows immediately. The variation of constants formula can be verified by direct derivation. \qed

The result below can be found  in \cite{Rall}.

\begin{proposition}\label{Prop-Rall}
Let $Y$ and $Z$ be normed linear spaces, $F:Y\to Z$ a map and suppose that the Gateaux derivative of $F$, $DF:Y\to\mathcal{L}(Y,Z)$ exists and is continuous at $y\in Y$. Then the Frech\`et derivative $F' $ of $F$ exists and is continuous at $y$.
\end{proposition}

\begin{proposition} \label{C1flow}
 Suppose, in addition to the hypotheses of Corollary \ref{globalexist}
  that the function $g$ is  continuously differentiable on $\mathbb{R} \times
 \R$. Then, for each $t\in\mathbb{R}$ the function $F(t,\cdot)$ is continuously Frech\'et differentiable on
$X $ with derivative given by
\begin{eqnarray*}
DF(t,u)v(x):=
\begin{cases}
 D_2 g(t, Ku(t,x)) Kv(t,x), & x\in \Omega,\\
 0, & x\in \mathbb{R}^N\backslash\Omega.
 \end{cases}
\end{eqnarray*}
where $D_2 g(t,y) = \partial_y g$ stands for the partial derivative of $g$ with respect to the second variable.
\end{proposition}
\proof From a simple computation, using the fact $g$ is continuously differentiable on $\mathbb{R}$, it follows that the Gateaux's derivative of $F(t,\cdot)$ is given by
\begin{eqnarray*}
DF(t,u)v(x):=
\begin{cases}
 D_2 g(t,  Ku(t,x))Kv(t,x), & x\in \Omega,\\
 0, & x\in \mathbb{R}^N\backslash\Omega.
 \end{cases}
\end{eqnarray*}

  The operator $DF(t,u)$ is clearly  a linear operator in $X$.  Using
 estimate (\ref{estimateLq}) and the continuity of $D_2 g$ it follows
  that,
 for any  $t \in \R$ and $u \in X$,  $ D_2 g(t, Ku(t,\cdot)) $ is bounded by a constant $k_3$.
  Using (\ref{estimateL1}),  we obtain
\begin{eqnarray*}
\| DF(t,u)v\|_{X} & \leq &
k_3 \|Kv\|_{L^p(\Omega)} \\
& \leq & k_3   \|v\|_{X},
\end{eqnarray*}
proving that  $DF(t, \cdot)$ is a bounded operator.

 We now prove that  $DF(t,\cdot): X \to \mathcal{L}(X)$ is  continuous.
If  $w$ and  $u$ in $X$, it follows from  (\ref{estimateLq}) that
\begin{eqnarray*}
& &\|DF(t,u  )v-DF(t,w  )v\|_{L^p(\Omega)}  \\
& \leq &\left\{ \int_{\Omega} |D_2 g(t, Ku  (x)) -  D_2 g(t, Kw
(x))|^p |Kv(x)|^{p}dx \right\}^{1/p}\\
& \leq &   \|J \|_{q}
\left\{ \int_{\Omega} |D_2 g(y, Ku  (x)) -  D_2 g(t, Kw  (x))|^p \, d \, x \right\}^{1/p} \| v\|_{L^p (\Omega) }.
 \end{eqnarray*}

  Now, from (\ref{estimateLq}),  it
 follows that   $Ku(x)  $  and $ Kw(x) $ are  in a bounded set of $\R^N$  for
 $w$ in a neighborhood of $u$ in $L^p(\Omega)$ and any  $t \in \R$, $x \in \Omega$ . Also,   $|   Ku (x) -  Kw  (x)| \to 0 $
 uniformly in   $ \Omega$,   as
 $ w \to u $ in $L^p(\Omega)$.
 Therefore,   $| D_2 g(t,  Ku  (x)) -  D_2 g(t,  Kw  (x))| \to 0 $
 uniformly in   $ \Omega$,   by the
 continuity of $D_2 g$.  From this the continuity of
  $DF(t,\cdot)$ at $u$ follows immediately.


Hence, it follows from Proposition \ref{Prop-Rall} that $F(t,\cdot)$ is Frechet differentiable with continuous derivative in $X$. \qed

\begin{remark}
Since the right-hand side of (\ref{CP}) is a $C^1$ function, the process generated by (\ref{CP}) in $X$  is  $C^1$  with respect to initial conditions.
\end{remark}

 From the results above, we have that, for each $\tau\in\mathbb{R}$ and $u(\tau,\cdot)\in X$, the unique solution of (\ref{CP}) with initial condition $u(\tau,\cdot)$ exists for all $t\geq\tau$ and   this solution $(t,x)\mapsto u(t,\tau,x;u_{\tau})$ (defined by (\ref{EP_1})) gives rise to a family of nonlinear $C^1$ evolution process on $X$ given by
\[
S(t,\tau)u_{\tau}(x):=u(t,\tau,x;u_{\tau}),\ t\geq\tau\in\mathbb{R},\ x\in\mathbb{R}^N.
\]

\section{Existence and regularity of the pullback attractor}\label{PullAttractors}

We   prove the existence of a pullback attractor $\{ \mathcal{A}(t); t\in\mathbb{R}\}$ in $X$ for the
evolution process $ \{S(t,\tau); t \geq \tau\in\mathbb{R}\} $ when $ 1 \leq p
 \leq  \infty $. In the estimates below the constant, $ |\Omega|^{\frac{1}{p} }$ should be taken as  $1$ when $p =\infty$.

\begin{lemma}\label{L_PullAbs}
Suppose that the hypotheses of Proposition \ref{C1flow} hold with the constant
 $k_1$ in (\ref{dissip2}) satisfying  $k_{1}< 1$
 Then the ball of $L^p(\Omega), \ 1 \leq p \leq  \infty$, centered
at the origin with  radius
$ (1+\delta)\frac{k_2 |\Omega|^{\frac{1}{p} }}{1-k_1}$,
 where $k_1$ and  $k_2$ are the constants appearing in (\ref{dissip2})
 and $\delta$ is any positive number
  pullback absorbs bounded subsets of $X$ under
 the  evolution process $\{S(t,\tau); t\geq\tau\in\mathbb{R}\}$ generated by (\ref{CP}) (with $|\Omega|^{\frac{1}{p} }$ replaced by $1$ if $p = \infty$).
\end{lemma}
 \proof  If  $u(t,\tau,x;u_{\tau})$ is a solution of (\ref{CP}) with
initial condition $u_{\tau}$ then, for $ 1 \leq p <\infty$
\begin{eqnarray*}
\frac{d}{dt}\int_{\Omega}|u(t,\tau,x;u_{\tau})|^{p}dx &=&\frac{d}{dt}\int_{\Omega}|u(t,\tau,x;u_{\tau})|^{p}dx\\
 & = & \int_{\Omega} p |u|^{p-1} \sgn(u)
   u_t(t,\tau,x;u_{\tau})   \,dx  \\
&=& -p\int_{\Omega}|u|^{p}(t,\tau,x;u_{\tau})dx +
p\int_{\Omega}|u|^{p-1} \sgn(u) g(t,Ku)dx.
\end{eqnarray*}

Using H\"{o}lder's inequality, condition (\ref{dissip1}) and  estimate
 (\ref{estimateLq}),  we obtain
\begin{eqnarray*}
\int_{\Omega}|u|^{p-1} \sgn(u) g(t,Ku)dx  &\leq&
\left(\|u(t,\tau,\cdot;u_{\tau})\|^{(p-1)q}_{L^{p}(\Omega)}\right)^{\frac{1}{q}}
 \left(
\int_{\Omega} | g(t,Ku)|^p dx\right)^{\frac{1}{p}}  \\
 &\leq&
\left(\|u(t,\tau,\cdot;u_{\tau})\|^{p}_{L^{p}(\Omega)}\right)^{\frac{p-1}{p}}
 \left(
\int_{\Omega} (k_1 |Ku| + k_2)^p dx\right)^{\frac{1}{p}}  \\
&\leq & \|u(t,\tau,\cdot;u_{\tau})\|^{p-1}_{L^{p}(\Omega)} \left[k_1
\|u(t,\tau,\cdot;u_{\tau})\|_{L^{p}(\Omega)}+k_2 |\Omega|^{\frac{1}{p}},\right]
\end{eqnarray*}
where $q$ is the conjugate exponent of $p$.

Hence
\begin{eqnarray*}
\frac{d}{dt}\|u(t,\tau,\cdot;u_{\tau})\|_{L^{p}(\Omega)}^{p} &\leq&
-p\|u(t,\tau,\cdot;u_{\tau})\|_{L^{p}(\Omega)}^{p}+pk_1
\|u(t,\tau,\cdot;u_{\tau})\|_{L^{p}(\Omega)}^{p}\\
& +&
pk_2{|\Omega|}^{\frac{1}{p}} \|u(t,\tau,\cdot;u_{\tau})\|^{p-1}_{L^{p}(\Omega)}\\
&=& p\|u(t,\tau,\cdot;u_{\tau})\|_{L^{p}(\Omega)}^{p} \left[-1 +
k_1  +
\frac{k_2 |\Omega|^{\frac{1}{p}}}{\|u(t,\tau,\cdot;u_{\tau})\|_{L^{p}(\Omega)}}
\right].
\end{eqnarray*}
Since $k_1  <1$,  $\varepsilon =1-k_1
 >0$. Thus, while
$\|u(t,\tau,\cdot;u_{\tau})\|_{L^{p}(\Omega)} \geq
 (1+\delta) \frac{k_2|\Omega|^{\frac{1}{p}}}{\varepsilon}$, we have
\begin{eqnarray*}
\frac{d}{dt}\|u(t,\tau,\cdot;u_{\tau})\|_{L^{p}(\Omega)}^{p} &\leq&
p \|u(t,\tau,\cdot;u_{\tau})\|_{L^{p}(\Omega)}^{p}(-\varepsilon +
\frac{\varepsilon}{1+\delta})\\
&=&- p\frac{\delta}{1+\delta}\varepsilon \|u(t,\tau,\cdot;u_{\tau})\|_{L^{p}(\Omega)}^{p}.
\end{eqnarray*}
Therefore,  while
$\|u(t,\tau,\cdot;u_{\tau})\|_{L^{p}(\Omega)} \geq
(1+\delta)\frac{k_2|\Omega|^{\frac{1}{p}}}{1-k_1 }$, we have
\begin{eqnarray} \label{boundsol}
\|u(t,\tau,\cdot;u_{\tau})\|_{L^{p}(\Omega)}^{p} &\leq&
e^{-\frac{\varepsilon \delta p}{(1+\delta)}t} \| u_{\tau}\|_{L^{p}(\Omega)}^{p} \nonumber \\
&=&e^{- \frac{\delta p}{(1+\delta)}(1-k_1)
t}\| u_{\tau}\|_{L^{p}(\Omega)}^{p}.
\end{eqnarray}
 From this, the  result follows easily for $1 \leq p \leq \infty$. Since the estimates are uniform in $p$, it also follows for $p = \infty$, by taking the limit with $p \to \infty$.  \qed

\begin{theorem}\label{Theor1}
In addition to the  conditions of  Lemma \ref{L_PullAbs}, suppose that
 $ \| J_x(x, \cdot) \|_{L^q(\Omega)} $ is bounded, where $q$ is the conjugate exponent of $p$ and $D_2g(t,\cdot)$ the derivative of $g$  with respect to the
 second variable is bounded in bounded sets of $\mathbb{R}$, uniformly for $t \in \mathbb{R}$.  Then there exists a pullback attractor $\{ \mathcal{A}(t); t\in\mathbb{R}\}$ for the evolution process $\{S(t,\tau); t\geq\tau\in\mathbb{R}\}$ generated by (\ref{CP}) in $X= L^p(\Omega)$  and
  $\mathcal{A} (t )$ is contained in the ball of radius
  $ \frac{k_2|\Omega|^{\frac{1}{p}}}{1-k_1 }$ in $L^p(\Omega)$,  for any $t \in \mathbb{R}$ and $1 \leq p \leq \infty$.
\end{theorem}
\proof
 From  (\ref{EP_1}), it follows that
\begin{equation*}
S(t,\tau)u_{\tau} =T(t,\tau) u_{\tau}  + U(t,\tau) u_{\tau}
 ,\ t\geq\tau,\ x\in\Omega
\end{equation*}
where
\[
T(t,\tau)u_{\tau}(x):=e^{-(t-\tau)} u_{\tau}(x)
\]
 and
\[
U(t,\tau)u_{\tau}(x):=\int_{\tau}^{t}e^{-(t-s)} g(s, Ku(s,\tau, x;u_{\tau}))\, ds.
\]

Suppose $ u_{\tau} \in B$, where $B$ is a bounded subset of
$X$. We may suppose that $B$ is contained in the ball  centered at the origin of radius $r>0$. Then
\[
\|T(t,\tau) u_{\tau}\|_{L^p(\Omega)}\leq r e^{-(t-\tau)},\ t\geq\tau\in\mathbb{R}.
\]

Also, from (\ref{boundsol}), we have that $\|u(t,\tau,\cdot;u_{\tau})\|_{L^p(\Omega)}\leq M$, for $t\geq \tau$, where $M=\max \{ r ,\frac{2k_2|\Omega|^{\frac{1}{p}}}{1-k_1}$\}.  Hence, for all $t \geq \tau$,
 and $x \in \Omega$,  we obtain
\begin{equation*}
|\partial_x U(t,\tau)u(\tau,x)|\leq\int_{\tau}^{t}e^{-(t-s)}|D_2
g(s, Ku(s,\tau, x;u_{\tau}))||\partial_xKu(s,\tau,x;u_{\tau})|\, ds.
\end{equation*}
{   By H\"{o}lder's inequality we
have that
\begin{equation}\label{Est-No}
\begin{split}
|\partial_xKu(s,\tau,x;u_{\tau})|&   \leq\int_{\Omega}
  |\partial_x J(x,y) |   |u(s,\tau, y ;u_{\tau}) | \, d \, y \\
 &  \leq
   \|\partial_x J(x,\cdot)\|_{L^q(\Omega)}\|u(s,\tau,\cdot;u_{\tau})\|_{L^p(\Omega)} \\
  &\leq C \|u(s,\tau,\cdot;u_{\tau})\|_{L^p(\Omega)}
\end{split}
\end{equation}
for some $C>0$ such that $\sup_{x\in\Omega}\|\partial_x J(x,\cdot)\|_{L^q(\Omega)}\leq C$.  Since  $\|u(t,\tau,\cdot;u_{\tau})\|_{L^p(\Omega)}\leq M$, for $t\geq \tau$,
 we obtain,  using
 estimate (\ref{estimateLq}) that  $ Ku(t,\tau,\cdot;u_{\tau})$ is bounded by a constant independent of $ t \geq \tau $ and $u_{\tau} \in B$.  Using the hypotheses, it follows
  that $ D_2 g(t,  Ku(t,\tau,\cdot;u_{\tau})) $ is bounded by a constant $k_3$  for any   $t \geq \tau$  and $u_{\tau} \in B$.   Thus
\[
\begin{split}
|\partial_x U(t,\tau)u_{\tau}(x)|&\leq k_3 C  M \int_{\tau}^{t}e^{-(t-s)} \, ds\\
&\leq  k_3 C  M.
\end{split}
\]
}

Therefore, for $t\geq\tau$ and any $ u_{\tau} \in B$, the value of
$\|\partial_x U(t,\tau)u_{\tau} \|_{L^{p}(\Omega)}$ is bounded by
a constant (independent of $\tau$ and $u_{\tau} \in B$). It follows that  for all $u_{\tau} \in B,$
we have that $U(t,\tau) u_{\tau} $ belongs to  a  (fixed)
ball of $W^{1,p}(\Omega)$  for all $u_{\tau} \in B$.
 From Sobolev's Embedding Theorem, it  follows that
  $U(t,\tau) $ is a compact operator, for any $ \tau \leq t$ and then, from
Theorem \ref{critcompact}, it follows that the process  $\{S(t,\tau); t\geq\tau\in\mathbb{R} \} $ is
 asymptotically compact.


Therefore, it follows from Lemma \ref{L_PullAbs} and
 Theorem \ref{Teorema 0.2.1} that the
pullback attractor $\{ \mathcal{A}(t); t\in\mathbb{R}\}$ exists and
is the pullback $\Omega$-limit set of any bounded subset of $X$
containing $\mathcal{B_{\delta}}$, where
$\mathcal{B}_{\delta}=
B\left(0, \frac{ (1+ \delta)k_2|\Omega|^{\frac{1}{p}}}{1-k_1}\right)$,
 for any $\delta > 0$ i.e., for instance
\begin{equation*}
\mathcal{A}(t)=\Omega_\wp(\mathcal{B}_{\delta},t):=\bigcap_{s\leq
t}\overline{\bigcup_{\tau\leq s}S(t,\tau)\mathcal{B}_{\delta}},\ \ \mbox{for
each}\  t\in\mathbb{R}.
\end{equation*}
From  this, since $\mathcal{B}_{\delta} $  pullback absorbs bounded subsets of   $X$,  it also follows that
 $ \mathcal{A}(t)$  is contained in the ball centered at the origin of radius $\frac{k_2|\Omega|^{\frac{1}{p}}}{1-k_1\|\beta\|_{L^{\infty}(\mathbb{R})}}$ in $L^{p}(\mathbb{R}^{N})$, for any $t \in \R$  $1 \leq p \leq \infty$.
\qed

\begin{theorem}
Assume the same conditions as in Lemma \ref{L_PullAbs}.
 Then there exists a bounded set of $W^{1,p}(\Omega)$, $1 \leq p \leq \infty$ containing
  $\mathcal{A}(t)$, for any $ t \in \R$.
\end{theorem}

\proof From  Theorem  \ref{Theor1}, we obtain that
 $\mathcal{A}(t)$
is contained in the ball centered at the origin and radius  $ \frac{k_2|\Omega|^{\frac{1}{p}}}{1-k_1 }$  in $L^{p}(\Omega)$. Furthermore, from Theorem \ref{Theor1}, we get that the pullback attractor can be written as the set of all global bounded solutions. Hence, if $u(t,\tau,x;u_{\tau})$ is a solution of (\ref{CP}) such that $u(t,\tau,x;u_{\tau})\in\mathcal{A}(t)$ for all $t\in\mathbb{R}$, then
\begin{equation}\label{equality_Lp}
u(t,\tau,x;u_{\tau})=\int_{-\infty}^{t}e^{-(t-s)}g(s, Ku(s,\tau,x;u_{\tau}))\, ds,
\end{equation}
where the equality above is in the sense of $L^{p}(\mathbb{R}^{N})$.

{
Proceeding as in the proof of  the Theorem \ref{Theor1} (see the estimate \eqref{Est-No} above), we have that
\[
\begin{split}
|\partial_{x} u(t,\tau,x;u_{\tau})| &\leq
\int_{-\infty}^{t}e^{-(t-s)}
|D_2 g(s, Ku(s,\tau,x;u_{\tau}))||\partial_{x}Ku(s,\tau,x;u_{\tau})|\, ds\\
&\leq k_{3} C M.
\end{split}
\]
for some $C>0$ such that $\sup_{x\in\Omega}\|\partial_x J(x,\cdot)\|_{L^q(\Omega)}\leq C$, and $k_3$ is a bound  of
  $ |D_2 g(t,  Ku(\cdot))| $   for  $u$
 in the ball of radius   $ \frac{2k_2|\Omega|^{\frac{1}{p}}}{1-k_1 }$ in $L^p(\Omega)$.  It follows that
$\mathcal{A}(t)=S(t,\tau)\mathcal{A}(\tau)$  is in a bounded set of
$W^{1,p}(\Omega)$  uniformly for  $t\in\mathbb{R}$.\qed }

\section{A comparison result} \label{comparison}

In this section we prove a comparison result in
$L^\infty(\mathbb{R}^N)$ and  use  it to prove  the invariance of
some sets  under  $S(t,\tau)$.


\begin{definition}
A function $v=v(t,\tau,x;v_\tau)$ is a sub solution of the Cauchy problem
(\ref{CP}) with initial condition $v_{\tau}$ in
$L^\infty(\mathbb{R}^N)$ if $v(\tau,x) = v_{\tau}(x)$ for
almost everywhere (\textrm{a.e.}) $x\in \Omega$, $v$ is continuously
differentiable with respect to $t$ and satisfies
\begin{equation*}
\partial_tv\leq -v+g(t, Kv)
\end{equation*}
a. e. in $[\tau,+\infty]\times\Omega$.
\end{definition}

Analogously, the function $V(t,\tau,x;u_\tau)$ is a super solution of (\ref{CP})
with initial condition $V_\tau$
if $V(\tau,x)= V_\tau(x)$, for
almost everywhere (\textrm{a.e.}) $x\in \Omega$, $V$ is continuously
differentiable with respect to $t$ and satisfies
\begin{equation*}
-V+g(t, K V)\leq\partial_t V.
\end{equation*}

\begin{theorem}\label{Comparation}
Assume the functions $f, g$ and $h$ satisfy  the
same conditions required for $g$ in  Proposition \ref{C1flow}. Assume also that
 $f(t,x) \leq g(t,x) \leq h(t,x)$, a.e.  for $(t,x) \in \R \times \R$ and
 the functions $f(t, \cdot)$ and $h(t, \cdot)$ are    increasing.
 Let $v(t,\tau,x;u_\tau)$, $(V(t,\tau,x;u_\tau))$ be a sub
solution (super solution) of the Cauchy problem \eqref{CP} with
$f$ ($h$) substituted for $g$ and
initial condition $v_\tau (V_\tau) \in L^\infty(\Omega)$ and
 $u(t,\tau,x;u_\tau)$ a solution of \eqref{CP}, with initial condition
$u_\tau  \in L^\infty(\Omega)$ . Then, if
$v(\tau, \cdot) \leq u(\tau, \cdot)  \leq V(\tau, \cdot) \, \mathrm{a. e.\ in \ }  \Omega. $
 $$
v(t,\tau,x;v_\tau)\leq u(t,\tau,x;u_\tau) \leq V(t,\tau,x;V_\tau) \
\mathrm{a. e.\ in}\ [\tau,+\infty]\times\Omega.
$$
\end{theorem}

\proof
 Consider the operator $G$ on
$L^{\infty}([\tau, T]\times\Omega)$, for some $T\in
\mathbb{R}$ (to be fixed later), defined by
\begin{eqnarray} \label{defG}
(G\phi)(t,x)=
\begin{cases}
e^{-(t-\tau)}\phi(\tau,x)+\displaystyle\int_{\tau}^{t}e^{-(t-s)}g(s, K\phi)(s,x)
)\, ds,
& x\in \Omega,\\
 0, & x\in \mathbb{R}^N\backslash\Omega.
 \end{cases}
\end{eqnarray}

If $T > \tau$, it follows  from (\ref{estimateLq}) that
$  | K\phi(t,x)| \leq
\|\phi\|_{L^{\infty}([\tau, T]\times\Omega) } $ and, therefore,
$|g(t, K\phi(t,x))  |$ is bounded by a constant $b$.

Thus, for any $(t,x) \in [\tau, T] \times \Omega $
\[
\begin{split}
|(G\phi)(t,x) |&\leq e^{-(t-\tau)}|\phi(\tau,x)|+\int_{\tau}^{t}e^{-(t-s)}|
g(s,K\phi(s,x)) |\, ds\\
&\leq e^{-(t-\tau)}|\phi(\tau,x)|+ b \int_{\tau}^{t}e^{-(t-s)} \, d\, s \\
& \leq \|\phi(\tau,\cdot) \|_{L^\infty(\Omega)} + b (T-\tau)
\end{split}
\]

 Therefore  $G$ is well defined as an operator on
$L^{\infty}([\tau,T]\times\Omega)$ and
furthermore $(G\phi)(\tau,x) =\phi(\tau,x)$.

Let
 $u_{\tau} \in L^\infty(\Omega)$, $L \in \R$ a positive number, and consider
 the set $\mathfrak{M_{u_{\tau} }} $
in $L^{\infty}([\tau,T]\times\Omega)$ given by

 \[ \mathfrak{M}_{u_{\tau} }
= \{\phi; \  \phi(\tau, x) = u_{\tau}(x) \textrm{ and }
 |\phi(t, x) - u_{\tau}(x) | \leq L, \textrm{ for }  \tau \leq t \leq T  \}.
 \]
if $\phi \in \mathfrak{M}_{u_{\tau}}$ then,
 from  (\ref{estimateLq}),  it follows that, for all
 $ (t,x) \in  [\tau,T]\times\Omega $ and $\phi \in \mathfrak{M}_{u_{\tau}}$
\begin{equation} \label{estimateKphi}
 | K\phi(t,x)| \leq
 \|\phi\|_{L^{\infty}([\tau, T]\times\Omega) }
 \leq   (L+ \|u_{\tau}\|_{L^{\infty}(\Omega)})
 \end{equation}
 and, therefore,
$|g(t, K\phi(t,x))|$ is bounded by a constant $b$ for
  $ (t,x) \in  [\tau,T]\times\Omega $ and $\phi \in \mathfrak{M}_{u_{\tau}}$.
Thus, for any $(t,x) \in [\tau, T] \times \Omega $
\[
\begin{split}
|(G\phi)(t,x) -u_{\tau}(x)|  & \leq (1-e^{-(t-\tau)})| u_{\tau}(x)|+
\int_{\tau}^{t}e^{-(t-s)}|g(s, K\phi(s,x)) |\, ds \\
&\leq (1-e^{-(t-\tau)})| u_{\tau}(x)|+ b \int_{\tau}^{t}e^{-(t-s)} \, d\, s \\
& \leq  (1-e^{-(t-\tau)}) \|u_{\tau}\|_{L^{\infty}(\Omega)}  + b (T-\tau).
\end{split}
\]
and $(G\phi)(t,x)$ belongs to $\mathfrak{M}_{u_{\tau}}$, for any
 $\phi \in  \mathfrak{M}_{u_{\tau}}$ if $(T-\tau)$ is small enough.

 We now prove that $G$ is a contraction in $\mathfrak{M}_{u_{\tau}}$
for some $T> \tau$.
Using \eqref{estimateKphi} and the fact that
 $g$ is locally Lipschitz in the second variable, we obtain the existence of
 a constant
 $k_{\mathfrak{M}}$, such that, for $\phi_1, \phi_2 \in \mathfrak{M}_{u_{\tau}}$
\begin{eqnarray*}
|g(t, K\phi_1(t,x)) -  g( t, K\phi_2(t,x))|
& \leq &  k_{\mathfrak{M}}  |K\phi_1(t,x)   -   K\phi_2(t,x)| \\
& \leq &  k_{\mathfrak{M}}  |\phi_1(t,x)   -  \phi_2(t,x)|.
 \end{eqnarray*}

Therefore
\[
\begin{split}
|(G\phi_1)(t,x) & -(G\phi_2)(t,x)  |\\
&\leq \int_{\tau}^{t}e^{-(t-s)}|g(t, K\phi_1(s,x) )
-g(s, K\phi_2(s,x))|\, ds\\
&\leq k_{\mathfrak{M}}
\int_{\tau}^{t}e^{-(t-s)} |\phi_1(s,x)   -   \phi_2(s,x)|  \, ds\\
&\leq k_{\mathfrak{M}}
 \|\phi_1   -  \phi_2\|_{L^{\infty}([\tau, T]\times\Omega)  }
\int_{\tau}^{t}   e^{-(t-s)}   \, ds \\
&\leq k_{\mathfrak{M}}  (T-\tau)
 \|\phi_1   -  \phi_2\|_{L^{\infty}([\tau, T]\times\Omega)  }. \\
\end{split}
\]

Hence $G$ is a contraction in $\mathfrak{M}_{u_{\tau}}$ for $T-\tau$ small enough.

 If $\phi(t,x)$ belongs to   $\mathfrak{M}_{u_{\tau}}$ and
 $z(t,x)=\lim_{n\to \infty}(G^{n} \phi)(t,x)
$
then,  from the continuity of $G$, we obtain that
$$
(Gz)(t,x)=G\left(\lim_{n\to \infty} (G^{n} \phi)(t,x)\right)=\lim_{n\to
\infty} (G^{n+1}\phi)(t,x)=z(t,x).
$$

Therefore  $z(t,x)$ is (the unique)   fixed point of $G$ in
 $\mathfrak{M}_{u_{\tau}}$ and it follows that  $ z(t,x) =   u(t,\tau,x;u_\tau)$,
 the unique  solution of (\ref{CP}) in
$[\tau, T]\times\Omega$ with initial condition
$z(\tau,x)=  u_{\tau}(x)$. In particular
\[
u(t,\tau,x;u_\tau)=\lim_{n\to \infty} (G^{n}u_{\tau})(t,x)\ \ \mathrm{in}\ \
L^{\infty}([\tau,T]\times\Omega),
\]
where, by an abuse of notation, we still denote by $u_{\tau}$  the function in $\mathfrak{M}_{u_{\tau}}$, which is equal to $u_{\tau}$ for $ \tau \leq t \leq T$.

 Let $F$ and $H$ defined as in \eqref{defG}   with $g$ replaced
 by $f$ or $h$, respectively.
The same results proved above for $G$ remain true for $F$ and $H$ and
 furthermore, since $g$ and $h$ are increasing  it follows that
    $F$ and $G$  are now  monotonic, that is,  for
any $\phi_{1},\phi_{2} \in L^{\infty}([\tau,T]\times\Omega)$
with $\phi_{1} \leq \phi_{2}$ a.e. in
$[\tau,T]\times\Omega$, we have  $F(\phi_{1})\leq F(\phi_{2})$ and
$H(\phi_{1})\leq H(\phi_{2})$
a.e. in $[\tau,T]\times\Omega$ (for some $T>\tau$).

 Now, if $v$ is a sub solution of (\ref{CP}), with $g$ replaced by $f$
 and   initial condition
 $v_\tau \leq u_\tau  $ a.e.,
  we obtain
$$
v(t,\tau,x;v_\tau)\leq
e^{-(t-\tau)}v(\tau,x)+\int_{\tau}^{t}e^{-(t-s)}f(s, K v(s,\tau,x;v_\tau))\, ds\,\,
\textrm{a.e. in}\ [\tau,T]\times\Omega.
$$

Hence
\[
v(t,\tau,x;v_\tau)\leq (Fv)(t,\tau,x;u_\tau)\,\,
\textrm{a.e. in}\ [\tau,T]\times\Omega.
\]
and  it follows
that $v(t,\tau,x;v_\tau)\leq (F^{n}v)(t,\tau,x;v_\tau)$ a.e. in $[\tau,T]\times\Omega$. Thus, passing to the limit as above, with an eventually
 smaller $T$ to
 ensure that {  $v(t,\tau,x;v_\tau)$ belongs to  $\mathfrak{M}_{v_{\tau}}$ and }
  $H$ is a contraction in  $\mathfrak{M}_{v_{\tau}}$
 we  obtain
  \begin{equation}  \label{comparesub}
v(t,\tau,x;v_\tau)\leq
w(t,\tau,x;v_\tau),    \textrm{ \ a.e. in \  }  [\tau,T]\times\Omega
\end{equation}
 where    $w(t,\tau,x;v_\tau)$ is
 the unique  solution of (\ref{CP}) in
$[\tau, T]\times\Omega$ with initial condition
$   v_{\tau}$, and $f$ substituted for $g$  .

 Since  $ v(\tau,x) \leq u(\tau,x)$  $a.e.$ in $\Omega$,
 and $f \leq g$  it follows that

\[
(F^{n} {v_\tau})(t,x)\leq (G^{n}u_\tau)(t,x) \,\,
\textrm{a.e. in}\ [\tau,T]\times\Omega.
\]
Passing to the limit as $n$ goes to infinity, we obtain
 \begin{equation}  \label{comparesol}
w(t,\tau,x;v_\tau)\leq
u(t,\tau,x;u_\tau),   \textrm { \ a.e. in \  }  [\tau,T]\times\Omega.
\end{equation}

 From \eqref{comparesub} and \eqref{comparesol} it follows that
  \[
 v(t,\tau,x;v_\tau)\leq
u(t,\tau,x;u_\tau).
\]

If  $V(t,\tau,x;u_\tau)$ is a super solution of the Cauchy problem of (\ref{CP})  with initial condition $V_\tau$,  $ u_\tau \leq V_\tau $ we obtain, by analogous arguments
 \[
 u(t,\tau,x;u_\tau)\leq
 V(t,\tau,x;V_\tau),   \textrm { \ a.e. in \  }  [\tau,T]\times\Omega.
\]

 Since the estimates above are uniform in a bounded set containing
 the initial conditions, we
may  extend the result to the interval $[T, (2T-\tau)]$  and, by
iteration, we can complete the proof of the theorem.
\qed

 From Theorem \ref{Comparation}, we immediately obtain invariance of some subsets
 of $L^{\infty}(\Omega)$ (or the functions in $L^{\infty}(\mathbb{R}^N)$ vanishing outside $\Omega$,
 with the obvious identification) under the process   $S(t, \tau)$.

 \begin{corollary}
 Suppose $f(t,x) = f(x)$ and $h(t,x) = h(x)$ satisfy the hypotheses of Theorem  \ref{Comparation}
 and are independent of $t$. If
  $v$ ($V$)  is an  non negative (non positive)
 equilibrium  of \eqref{CP} with $f$ ($h$)
 substituted for $g$ then
the set \[ U=\{u\in  L^{\infty}(\Omega);\ v \leq u \leq V \}
\quad ( U=\{u\in  C^{0}(\mathbb{R}^N); v \leq u \leq 0 \} )\]
is positively
invariant under the process $S(t, \tau)$, that is $S(t,\tau)U\subset
U$ for all $t\geq \tau$.
\end{corollary}

\section{Upper semi-continuity of the pullback attractors}\label{Sec6}

{
A natural question  is the dependence of the pullback attractors
of (\ref{CP})  upon parameters. Suppose, for instance, that we have a family
 of functions $g_{\beta}(t,x)$, with $\beta$ in some normed space,
  $g_{\beta}(t,x) \to  g_{\beta_0}(t,x)$  as  $\beta \to \beta_0$ and  denote by
 $S_\beta(t,\tau)u_\beta(\tau,x) = u_\beta(t,\tau,x;u_{\tau}) $ the evolution process associated with the problem (\ref{CP}) with
 $g$ replaced by $g_\beta$. Under suitable hypotheses, we  prove
 the upper semi-continuity of the pullback attractors  as $\beta\to\beta_0$, for each $t \in \R$
that is , we show that
\[
\mathrm{dist}(\mathcal{A}_\beta(t),\mathcal{A}_{\beta_0}(t))\to0\quad\mbox{as}\quad \beta\to\beta_0\ \mbox{in}\ L^\infty(\mathbb{R})
\]
where $\{A_\beta(t);\ t \in \R \}$ denotes the pullback attractor of $S_\beta(t,\tau)$ in $X$, and  $\mathrm{dist}(\mathcal{A}_\beta(t),\mathcal{A}_{\beta_0}(t))$ is the Hausdorff semi-distance in $L^p(\Omega)$  defined by (\ref{DD}), for each $t\in\mathbb{R}$.

\begin{theorem}\label{TheoUpper01}
Assume the same  hypotheses of  Lemma \ref{L_PullAbs}, with $g$ replaced by
 $g_{\beta}$ for  $\beta$  in a neighborhood of $\beta_0$.  Suppose additionally
 that  the Lipschitz constants of the  functions $g_{\beta}$
 in bounded sets of $\R \times \R$ can be chosen independently of $\beta$ and
 $g_{\beta}(t,x) \to  g_{\beta_0}(t,x)$  as  $\beta \to \beta_0$ uniformly
 in bounded sets of $\R \times \R$.
  Then  $ S_\beta(t,\tau)(u_\tau)  \to   S_{\beta_0} (t,\tau)(u_{\tau})$
 \ in the
  $L^p(\Omega)-$norm, $1\leq p \leq \infty$,  uniformly for $t\in[\tau,c]$ with
  $c<\infty$ and $u_\tau$ in a  bounded subset of $X$.
\end{theorem}

\proof
Let $u_\beta(t,\tau,x;u_{\tau})=S_\beta(t,\tau)u_\beta(\tau,x) $ and $u_{\beta_0}(t,\tau,x;u_{\tau})=S_{\beta_0}(t,\tau)u_\beta(\tau,x) $ be the solutions of (\ref{CP})  with initial condition $u_\tau$, and parameters  $\beta$ and $\beta_0$, respectively, Then
\begin{eqnarray*}
&& (u_\beta -u_{\beta_0})(t,\tau,x;u_{\tau})\\
 &=&\int_\tau^t e^{-(t-s)}[g_{\beta}(s,Ku_\beta(s,\tau,x;u_{\tau}))-
g_{\beta_0}(s,Ku_{\beta_0})(s,\tau,x;u_{\tau}))]\, ds.
\end{eqnarray*}
 Summing and subtracting the term $\int_\tau^t e^{-(t-s)}g_{\beta}(s,Ku_{\beta_0}(s,\tau,x;u_{\tau})) \, ds$,  we get
\begin{eqnarray*}
&&  (u_\beta -u_{\beta_0})(t,\tau,x;u_{\tau})\\
 &=&\int_\tau^t e^{-(t-s)}[g_{\beta}(s,Ku_\beta(s,\tau,x;u_{\tau}))-g_{\beta}(s,Ku_{\beta_0}(s,\tau,x;u_{\tau}))]\, ds\\
 &+&\int_\tau^t e^{-(t-s)}[g_{\beta}(s,Ku_{\beta_0}(s,\tau,x;u_{\tau}))-
g_{\beta_0}(s,Ku_{\beta_0}(s,\tau,x;u_{\tau}))]\, ds
\end{eqnarray*}
for any $t\geq\tau$ and $x\in\Omega$.

 From estimates \eqref{boundsol} and  \eqref{estimateLq} it follows
 that  $ Ku_\beta(s,\tau,x;u_{\tau})$
 and   $ Ku_{\beta_0}(s,\tau,x;u_{\tau})$ remain in a bounded set
 $B \subset \R$ for any $ \tau\leq s \leq t$ and $x \in \Omega$. If $k_B$ denotes the Lipschitz constant of $g_{\beta}$ in $[\tau, t] \times B$ and
\[
\|g_{\beta}-g_{\beta_0}\|_{\infty}:= \sup \{ |g_{\beta}(s,x)-
g_{\beta_0}(s,x); \ \ (s,x) \in  [\tau, t] \times B  \} ,
\]
we obtain, using  \eqref{estimateL1}

\begin{eqnarray*}
&& \|(u_\beta -u_{\beta_0})(t,\tau,\cdot;u_{\tau})\|_X\\
& \leq &
\int_\tau^t e^{-(t-s)}\|g_{\beta}(s,Ku_\beta(s,\tau,\cdot;u_{\tau}))-g_{\beta}(s,Ku_{\beta_0}(s,\tau,\cdot;u_{\tau})) \|\, ds\\
 &+&\int_\tau^t e^{-(t-s)} \|g_{\beta}(s,Ku_{\beta_0}(s,\tau,\cdot;u_{\tau}))-
g_{\beta_0}(s,Ku_{\beta_0}(s,\tau,\cdot;u_{\tau}))\|\, ds \\
& \leq &
k_B\int_\tau^t e^{-(t-s)}\| Ku_\beta(s,\tau,\cdot;u_{\tau})-Ku_{\beta_0}(s,\tau,\cdot;u_{\tau})
 \|\, ds\\
 &+&  \|g_{\beta}-g_{\beta_0}\|_{\infty} |\Omega|^{\frac{1}{p}} \int_\tau^t e^{-(t-s)} \, ds \\
& \leq &
k_B  \int_\tau^t e^{-(t-s)}\| u_\beta(s,\tau,\cdot;u_{\tau})-u_{\beta_0}(s,\tau,\cdot;u_{\tau})
 \|\, ds\\
 &+&  \|g_{\beta}-g_{\beta_0}\|_{\infty}  |\Omega|^{\frac{1}{p}}
\end{eqnarray*}
where we have indicated by $\| \cdot \|$ the norm in the space $X= L^p(\Omega)$.

By Gronwall's inequality, it follows that
\[
\|(u_\beta -u_{\beta_0})(t,\tau,\cdot;u_\tau)\|_{X}\leq
  |\Omega|^{\frac{1}{p}}  \| g_{\beta} -  g_{\beta_0}\|_{\infty}  e^{k_B (t-\tau)}
\]
concluding the proof.
\qed

\begin{theorem}
 Under  the same conditions of Theorem \ref{TheoUpper01},
the family of pullback attractors $\{\mathcal{A}_\beta(t); t\in\mathbb{R}\}$ is upper semicontinuous at $\beta=\beta_0$, at each $t \in \R$.
\end{theorem}
\proof
 Denote by $M $
 the radius of an  absorbing ball for the  process $S_{\beta}$
 for $\beta$ close to $\beta_0$,   whose existence is guaranteed by
 Lemma \ref{L_PullAbs}.
Given  $\delta > 0$, let $\tau\in\mathbb{R}$ such that $\mathrm{dist}(S_{\beta_0}(t,\tau)B\left(0;M \right),\mathcal{A}_{\beta_0}(t))<\delta/2$.
 By Theorem \ref{TheoUpper01}
$$
\|(u_\beta -u_{\beta_0})(t,\tau,\cdot;u_\tau)\|_{X}\ \to \ 0
$$
as $\beta\to\beta_0$  uniformly for $u_\tau$ in bounded subsets of $X$.
 Since  $\bigcup_{s\in\mathbb{R}}\mathcal{A}_{\beta} (s)\subset B(0;M)$, for
 $\beta $ in a neighborhood of $\beta_0$ , there exists $r_0>0$ such that
\[
\sup_{a_\beta\in\mathcal{A}_\beta(\tau)}\|S_\beta(t,\tau)a_\beta-S_{\beta_0}(t,\tau)a_\beta\|_{X}<\delta/2
\]
for all $\beta\in B(\beta_0,r_0)$, the ball centered at $\beta_0$ of radius $r_0>0$ in $L^\infty(\mathbb{R})$.

Then
\begin{eqnarray*}
&&\mathrm{dist}(\mathcal{A}_\beta(t),\mathcal{A}_{\beta_0}(t))\\
&\leq& \mathrm{dist}(S_\beta(t,\tau)\mathcal{A}_\beta(\tau),S_{\beta_0}(t,\tau)\mathcal{A}_\beta(\tau))+\mathrm{dist}(S_{\beta_0}(t,\tau)\mathcal{A}_\beta(\tau),S_{\beta_0}(t,\tau)\mathcal{A}_{\beta_0}(\tau))\\
&=&\sup_{a_\beta\in\mathcal{A}_\beta(\tau)} \mathrm{dist}(S_\beta(t,\tau)a_\beta,S_{\beta_0}(t,\tau)a_\beta)+\mathrm{dist}(S_{\beta_0}(t,\tau)\mathcal{A}_\beta(\tau),\mathcal{A}_{\beta_0}(t))\\
&<&\frac{\delta}{2}+\frac{\delta}{2}
\end{eqnarray*}
and the upper semicontinuity is proved. \qed }

 {\color{red}}
\section{The forward asymptotically autonomous case }\label{NS}

 We now consider the forward dynamics of (\ref{CP}) in the
 case where $g(t,x)$  converges to a function
 $g_0(x)$ independent of $t$, when $t \to \infty$.
  Then (with appropriate hypotheses),  the process generated by
(\ref{CP}) becomes an  \emph{asymptotically autonomous} process,
 which is defined as follows (see, for example \cite{mischaikow}).
\begin{definition}
 An evolution process  $\{ S(t,\tau); t \geq \tau\in\mathbb{R} \}$
 in a Banach space $X$ is called
 \emph{asymptotically autonomous}  with \emph{limit (autonomous) process}
 $\{S_0(t); t \geq 0 \}$
if
\[ S(t_j +\tau_j, \tau_j) x_j \to  S_0(t)x, \]
 for any three sequences $t_j \to t$,  $\tau_j \to \infty$, $x_j \to x$, as
$ j \to \infty$, with $x, x_j \in X, 0 \leq t, t_j < \infty $, and
 $\tau_j \geq t_0$.
\end{definition}

Consider the (autonomous) problem in $X = L^p (\Omega), 1 \leq p \leq \infty$.
\begin{equation}\label{probaut}
\begin{cases}
\partial_t u (t,x) =- u (t,x)  + g_0( Ku(t,x))  \ \ \textrm {for }
 \ \ t \geq  \tau \in \R \ \ \textrm{and}    \ \ x \in  \Omega, \\
u(\tau, \cdot) = u_\tau (\cdot)   \ \  \textrm{in}  \ \ \Omega, \\
u(t, x)   = 0  \ \  \textrm{for}   \ \ t \geq  \tau \in \R
 \ \ \textrm{and}  \ \  x \in
\mathbb{R}^N\backslash\Omega.
\end{cases}
\end{equation}

 \begin{lemma}
 Suppose $g(t,x)$ and $g_0$ satisfy the hypotheses required for $g$ in Lemma
 \ref{L_PullAbs} and   $ g(t,x) \to g_0(x) $ as $t \to \infty$, uniformly for
$x $  in bounded sets of $\mathbb{R}$.  Then the process   $ \Sprocess $ generated
 by (\ref{CP}) is asymptotically autonomous with limit autonomous process
 $\Sautprocess$ generated by
 (\ref{probaut}).
\end{lemma}
\proof
 Suppose $t_j \to t$,  $\tau_j \to \infty$, $u_j \to u$ in $X$, as $ j \to \infty$, with $u, u_j \in X$, $t\geq0$, and $t_j < \infty $.

 Then

 \begin{eqnarray} \label{asymptaut}
   S(t_j + \tau_j, \tau_j) u_j  -  S_0(t+ \tau_j, \tau_j) u &  =    &
     (S(t_j + \tau_j, \tau_j) u_j -  S(t + \tau_j, \tau_j) u_j) + \nonumber \\
    (S(t + \tau_j, \tau_j) u_j -  S(t + \tau_j, \tau_j) u)
  &  +  & (  S(t + \tau_j, \tau_j) u -  S_0(t+ \tau_j, \tau_j) u).
\end{eqnarray}

 If  $t_j \to t$,  the first term to the right  in
 \eqref{asymptaut} goes to  zero
 by continuity of the semigroup.
The second term also goes to zero, as $u_j \to u$ by continuity of the semigroup with respect to initial conditions. This can be proved in our case,
  as follows:    by Lemma \ref{L_PullAbs}, $S(s+\tau_j, \tau_j)u $  belongs to a bounded set in $L^p(\Omega)$, for
 $ s,  \tau_j \in [0, \infty)$ and $u$ a bounded set in $L^p(\Omega)$.
 Therefore, from estimate $(\ref{estimateLq})$ ,
 $K S(s+\tau_j, \tau_j)u_j$
  and $ KS_0(s+\tau_j, \tau_j)u$ belong to a bounded set in $L^{\infty} (\Omega)$,
  for
 $ s,  \tau_j \in [0, \infty$. Then, using  that $g$ is locally Lipschitz and
  Gronwall's inequality, we obtain that
 $ S(t + \tau_j, \tau_j) u_j  \to   S(t + \tau_j, \tau_j) u $ as claimed.

 For the last term, we may use  results of Section
 \ref{Sec6}.  Indeed, defining, for
 $\beta \in \R$, the function $g_{\beta}(t,x) = g(t+ \beta, x)$, it follows from the assumptions that $g_{\beta} (t,x) \to g_0(x)$ as $\beta \to \infty$, uniformly  in bounded sets of $\R \times \R$. Denoting by
  $ S_{\beta}(t , \tau)$ the semigroup generated by
 \eqref{prob}, with
 $g_{\beta}$ substituted for $g$, we obtain    from Theorem \ref{TheoUpper01},
 that  $ S_{\beta}(t , \tau)u  \to S_0(t, \tau)u$ as $\tau \to \infty$.
 Since   $ S_{\tau_{j}}(t , \tau)u = S(t + \tau_j, \tau)u $, the claimed
 convergence is proved.

 \qed

 \vspace{3mm}

 We now obtain some results concerning the convergence of solutions of
 (\ref{CP}) to equilibria of (\ref{probaut}).  Similar results were obtained in    (\cite{polacik}) and (\cite{mischaikow}),  but  we present here a more direct  approach for our specific problem, based on the properties of
an ``energy functional''  in
 $X$.

To define this functional we need   the following additional
 hypothesis on the function $g_0$:
\begin{itemize}
\item[\textit{(H1)}] There exists a constant $a>0$ such that
\[
|g_0(s)|\leq a\ \ \mathrm{for\ all}\ \ s\in\mathbb{R}.
\]
\item[\textit{(H2)}] The function $g_0$ is strictly increasing
 and the function
\[
f(s) = -\frac12 s^2 - i(s), \ \ s\in[-a,a],
\]
where $i$ is given by
\[
i(s) = -\int_{0}^{s} g_0^{-1}(\theta) d\theta, \quad s\in [-a,a],
\]
has a global minimum 
\end{itemize}

\vspace{3mm}

 We then define the functional in
$ Y := \{u \in X ;\  \  |u(x)| < a, \ \mbox{for all}\ x \in \Omega  \}$.

\begin{equation}\label{4.8}
\begin{split}
\mathcal{L}(u) & =  \int_{\R^n}[f(u(x))-f(\overline{u})] \, dx +
\frac{1}{4}\int_{\R^n}\int_{\R^n}J(x,y)[u(x)-u(y)]^{2} \, dxdy \\
 & =
  \int_{\Omega}[f(u(x))-f(\overline{u})] \, dx +
\frac{1}{4}\int_{\Omega}\int_{\Omega}J(x,y)[u(x)-u(y)]^{2} \, dxdy
\end{split}
\end{equation}
where  $\overline{u}$ is a global minimum
of $f$.

\begin{remark} If $|g(t,x)|\leq a$, for $(t,x) \in \R \times \R$
 the  hypotheses of   Lemma \ref{L_PullAbs}
 hold with $k_1 =0$ and $k_2 =a$.   It follows then,  from estimate
(\ref{boundsol}), that  the set $Y$ is (pullback and 
forward) absorbing and, therefore, 
  positively invariant for the process generated by (\ref{CP}).
 Also, since $g(t,x)  \to g_0(x)$ uniformly in bounded sets,  this property holds for  $x$  in a bounded set, if $t$ is sufficiently large.
 Therefore,   if we  are interested only   in solutions in $Y$, we may suppose that it holds for all  $x$ by doing  an appropriate ``cut-off'' of $g(t,\cdot)$, if necessary.



\end{remark}

\begin{proposition}
Under hypotheses (H1) and (H2), the
functional $\mathcal{L}$ is continuous (with respect to the
 $L^p $ norm).
\end{proposition}
\proof
Note that, given $u\in Y$, as $|u(x)|\leq a$, for all $x \in \Omega$ there exists a
positive constant $f_{0}$ such that
$$
|f(u(x))-f(\overline{u})|\leq |f(u(x))|+|f(\overline{u})|\leq f_{0},
\,\,\ \mbox{for almost every} \, (a.e.) \,\, x\in \Omega.
$$
Let $(u_{n})$ a sequence converging to $u$ in the norm of
$L^{p}(\Omega)$.Then, by Theorem 4.9 in \cite{Brezis}, we can extract   a subsequence $(u_{n_{k}})$, such that,
$u_{n_{k}}(x)\rightarrow u(x)$ $a.e.$ in $\Omega$. Since from
\textit{(H2)}, it follows that $f$ is continuous,
$f(u_{n_{k}}(x))\rightarrow f(u(x))$ $a.e$. Thus
$$
\lim_{k \rightarrow \infty}
[f(u_{n_{k}}(x))-f(\overline{u})]=[f(u(x))-f(\overline{u})], \,\,
a.e.
$$
and
$$
\lim_{k \rightarrow
\infty}[u_{n_{k}}(x)-u_{n_{k}}(y)]^{2}=[u(x)-u(y)]^{2}, \,\, a.e.
$$
Now, we write
$$
\mathcal{L}(u)=\mathcal{L}_{1}(u)+\mathcal{L}_{2}(u),
$$
where
$$
\mathcal{L}_{1}(u)= \int_{\Omega}[f(u(x))-f(\overline{u})]dx
$$
and
$$
\mathcal{L}_{2}(u)=\frac{1}{4}\int_{\Omega}\int_{\Omega}J(x,y)[u(x)-u(y)]^{2}dxdy.
$$
Since
\begin{eqnarray*}
|f(u_{n_{k}}(x))-f(\overline{u})| \leq f_{0} \in L^{1}(\Omega),
\end{eqnarray*}
by Lebesgue's dominated convergence theorem (see \cite{Brezis}), we have
\begin{equation}\label{L1}
\lim_{k \rightarrow \infty}
\mathcal{L}_{1}(u_{n_{k}})=\mathcal{L}_{1}(u).
\end{equation}
Analogously
$$
\left|u_{n_{k}}(x)-u_{n_{k}}(y)\right|^{2}\leq 4a^{2} \in
L^{1}(\Omega).
$$
Hence, also by Lebesgue's dominated convergence theorem, we obtain
\begin{equation}
\lim_{k \rightarrow \infty}
\mathcal{L}_{2}(u_{n_{k}})=\mathcal{L}_{2}(u).\label{L2}
\end{equation}
Therefore, from (\ref{L1}) and (\ref{L2}), it follows that
$$
\lim_{k \rightarrow \infty}
\mathcal{L}(u_{n_{k}})=\mathcal{L}(u).
$$
Thus  $(\mathcal{L}(u_{n}))$ is a sequence in $\mathbb{R}$ such that
every subsequence has a subsequence that converges to
$\mathcal{L}(u)$, and  we obtain
$$
\lim _{n\rightarrow \infty}\mathcal{L}(u_{n})=\mathcal{L}(u).
$$
  \qed

\begin{remark}
The integrand in the functional $\mathcal{L}$ given in (\ref{4.8}) is
always  non negative since $J$ is positive and
 $\overline{u}$ is a  global
minimum of $f$. Thus $\mathcal{L}$ is lower bounded.
\end{remark}

 As shown in \cite{FAH}   $\mathcal{L}$  is a Lyapunov functional for dynamical system
 $\Sdyn$ generated by (\ref{probaut}) in $ L^{2}(\Omega)$, with derivative along the solutions $u(t,\tau,\cdot;u_\tau)$ given by

$$
\frac{d}{dt}\mathcal{L}(u(t,\tau,\cdot;u_\tau)=-I(u(t,\tau,\cdot;u_\tau))
$$
where

$$
I(u(t,\tau,\cdot;u_\tau))=
\int_{\Omega}[Ku(t,\tau,x;u_\tau)-g_0^{-1}(u(t,\tau,x;u_\tau))]
[- u(t,\tau,x;u_\tau)
+ g_0(Ku(t,x;\tau,u_\tau))]dx.$$


 In the nonautonomous case, we can obtain a similar result, with an additional term which is small for $t$ big.
 The proof is also analogous to the one obtained in \cite{FAH} for the
 autonomous case and will be omitted.

\begin{theorem}
Suppose that the hypotheses (H1) and (H2) hold and the functions
 $g(t, \cdot)$ are invertible, for all $t \in \R$. Let
$u(t,\tau,\cdot;u_\tau)$ be a solution of (\ref{CP}) with
$u(t,\tau,x;u_\tau)\leq  a$. Then
$\mathcal{L}(u(t,\tau,\cdot;u_\tau))$ is differentiable with respect to
$t$ for $t > \tau$ and, for any $u \in L^{2}(\mathbb{R}^{N})$,
$$
\frac{d}{dt}\mathcal{L}(u(t,\tau,\cdot;u_\tau)=-I(u(t,\tau,\cdot;u_\tau)) +
R(t, u(t,\tau,\cdot;u_\tau))
$$
where
\[
\begin{split}
I(u(t,\tau,\cdot;u_\tau))&= \int_{\Omega}[Ku(t,\tau,x;u_\tau)-g_0^{-1}(u(t,\tau,x;u_\tau))]
[- u(t,\tau,x;u_\tau)\\
&+ g_0(Ku(t,x;\tau,u_\tau))]dx\\
R(t,u(t,\tau,\cdot;u_\tau))&=
\int_{\Omega}[g^{-1}(t, u(t,\tau,x;u_\tau)) -g_0^{-1}(u(t,\tau,x;u_\tau))]
[- u(t,\tau,x;u_\tau)\\
&+ g_0(Ku(t,\tau,x;u_\tau))]dx
\end{split}
\]
Furthermore, $I(u(t,\tau,\cdot;u_\tau))=0$ if and only if
$u(t,\tau,\cdot;u_\tau)$ is one of the equilibria of the
 autonomous system (\ref{probaut}).
\end{theorem}
\begin{remark}
If $u(t,\tau,\cdot;u_\tau)$ remains bounded, the second
 term $R(t, u(t,\tau,\cdot;u_\tau))$ in the expression above for the derivative
 of $\mathcal{L}(u(t,\tau,\cdot;u_\tau))$ goes to zero, as $t\to\infty$ and, therefore, the derivative along the solution of the nonautonomous system approaches $-I$, which is
 the derivative along a solution of the autonomous system.
\end{remark}

\vspace*{5mm}

 \begin{lemma}
  The forward orbits of the process $\Sprocess$ generated by (\ref{CP})
 are   precompact in $X$.
\end{lemma}
\proof
Follows from the existence of attractors, or may be proved using the variation of constants formula and the decomposition
 $S(t,\tau) u_{\tau} = T(t,\tau) u_{\tau} + U(t,\tau) u_{\tau} $ given
 in Theorem  \ref{Theor1}.
\qed

\begin{theorem} \label{lasalle} Suppose that $\mathcal{L}$ assumes a finite number of values in  the set $E$ of equilibria of
 the dynamical system (\ref{probaut}). If $  u(t,\tau,x ; u_{\tau}) $ is a solution of
  (\ref{CP}) then
  $$  u(t,\tau,x ; u_{\tau}) \to \mathcal{M}, \  as \    t \to  + \infty, $$
  where  $\mathcal{M}$ is a level set of the functional $\mathcal{L}$ restricted to  $E$.
\end{theorem}
\proof
Let
$\mathcal{O} =  S(t, \tau)u_{\tau}$ denote  the forward orbit
 of $u_{\tau}$   and
$E_1, E_2, \ldots, E_m$ the level sets of $\mathcal{L}$ in $E$, in ascending order.
 Observe that $\overline{\mathcal{O}} $ is  positively invariant and compact.

 Write $\tilde{E_j} = E_j \cap \bar{\mathcal{O} }, \quad i=1,2, \ldots, m$,
   $\tilde{E} = E \cap \bar{\mathcal{O} } = \cup_{j=1}^m \tilde{E_j}$.
 Then,  each   $\tilde{E_j}$ is compact, since $E_j$ is closed.
 We first claim that $\tilde{E}$ is not  empty.  In fact, if this were true,
 then $I  <0$ in the compact set  $\overline{\mathcal{O}} $
 implies  $I  < -m$, for some positive number $m$    and then $\frac{d}{dt} \mathcal{L}(  S(t, \tau) u_{\tau}) < -m/2$ in  $\bar{\mathcal{O}}$ for $t$ big enough and,  a fortiori, also in  ${\mathcal{O}}$. But this is impossible, since
$\mathcal{L}$
  is bounded below.

Let  then $\tilde{E}_k$ be the minimum level which intercepts
   $\overline{\mathcal{O}}$.    We  now show that no other level is in the $\Omega$-limit of $  \mathcal{O}$. For this,
  choose a  neighborhood $V_i$ of each  $\tilde{E}_i $, such that
 $ \max \{\mathcal{L}(u) \ : \ u \in V_i \} <
  \min \{\mathcal{L}(u) \ : \ u \in V_{i+1} \} $,  $i=1, 2, \ldots, m-1$.
Let  $\mathcal{L}_k$ be the value of $\mathcal{L}$ in $\tilde{E}_k$ and
$W^i_{\varepsilon} = \{u \in V_i \ :  \  \mathcal{L}_i - \varepsilon < \mathcal{L}(u) <  \mathcal{L}_i + \varepsilon
   \}$,  $i=1, 2, \ldots, m-1$, with $\varepsilon$ chosen in such a way that
  $\overline{W^i_{\varepsilon}} \subset V_i$.
 Reasoning again by compactness we may find   $\bar{t}$ such that
 $\frac{d}{dt} \mathcal{L}(  S(t, \tau) u_{\tau}) < 0$
 in   $\overline{\mathcal{O}} - \cup_i W^i_{\varepsilon}$,
 if $t \geq \bar{t}$.
 Now, if $s > \bar{t}$ is such that  $ S(s, \tau) u_{\tau}$ belongs to $W^k_{\varepsilon}$  (such an $s$ exists, since the $\Omega$-limit set intersects $E_k$), then
  $ S(t, \tau) u_{\tau}$ remains in $W^k_{\varepsilon}$  for all $t> s$.  In fact,
 $ \mathcal{L}(S(t, \tau) u_{\tau})$ cannot become bigger than $ \mathcal{L}_k + \varepsilon $,
 since it decreases in    $ V_i \setminus W^i_{\varepsilon}$. It cannot also become
 smaller than
  $ \mathcal{L}_k - \varepsilon $. In fact  in this case,
$ S(t, \tau) u_{\tau}$ will never enter
 any of the neighborhoods $\cup_i W_i$ again since  $ \mathcal{L}(S(t, \tau) u_{\tau})$ decreases in  $\overline{\mathcal{O}} - \cup_i W^i_{\varepsilon}$.
 This is a contradiction, since   it must return to $W^k_{\varepsilon}$.

  Therefore, we must have  $ S(t, \tau)u_{\tau} \subset  W^k_{\varepsilon}$
 for $t$ big enough. Since $ \varepsilon $ is arbitrary,  we must have
 $\cap_{t \geq \tau} \cup_{s \geq t} \overline{ S(s, \tau)( u_{\tau})} \subset E_k$, as claimed

\qed

\begin{corollary}
 Suppose, that (\ref{probaut}) has a finite number of equilibria in each level set of $\mathcal{L} $. Then
 any solution $u(t,\tau,x ; u_{\tau})$ of
  (\ref{CP}) converges to a single equilibrium of  (\ref{probaut}).
 \end{corollary}
 \proof
  We must show that the $\Omega$-limit of  $u(t,\tau,x ; u_{\tau})$  is
 a single equilibrium.
  From Theorem \ref{lasalle} it follows that it must be a subset of
 the set of $E$ equilibria of (\ref{probaut}) and also
 contained in  a level set of $\mathcal{L}$.
  Since it must also be connected, the result follows immediately.
 \qed
 
\end{document}